\documentclass[11pt]{article}

\usepackage{amsmath,amsfonts,amssymb,amsthm,bbm}
\usepackage{graphics,epsfig}
\usepackage{color}
\usepackage{array}
\usepackage{ulem}
\usepackage{mathrsfs}
\newtheorem{theorem}{Theorem}[section]
\newtheorem{proposition}{Proposition}[section]

\numberwithin{equation}{section}

\allowdisplaybreaks

\renewcommand\today{\number\year/\number\month/\number\day}

\def\qed{\hfill$\Box$\medskip}

\begin{document}

\noindent\makebox[60mm][1]{\tt {\large Version:~\today}}

\bigskip
\noindent{
{\Large\bf  Scaling limit of the local time of the Sinai's random
walk}\footnote{
\noindent  Supported by NSFC (NO.11131003) and  985 Program.
}
\\

\noindent{
Wenming Hong\footnote{ School of Mathematical Sciences
\& Laboratory of Mathematics and Complex Systems, Beijing Normal
University, Beijing 100875, P.R. China. Email: wmhong@bnu.edu.cn} ~  Hui Yang\footnote{ School of Mathematical Sciences
\& Laboratory of Mathematics and Complex Systems, Beijing Normal
University, Beijing 100875, P.R. China. Email: yanghui2011@mail.bnu.edu.cn} ~
Ke Zhou\footnote{ School of Mathematical Sciences
\& Laboratory of Mathematics and Complex Systems, Beijing Normal
University, Beijing 100875, P.R. China. Email:zhouke@mail.bnu.edu.cn}

\noindent{
(Beijing Normal University)
}

}

\vspace{0.1 true cm}

\begin{center}
\begin{minipage}[c]{12cm}
\begin{center}\textbf{Abstract}\end{center}
\bigskip
We prove that the local times of a sequence of Sinai's random walks   convergence to those    of Brox's diffusion by proper scaling, which is  accord with the result of Seignourel (2000). Our proof is based on   the convergence of the branching processes in random environment by Kurtz~(1979).

\mbox{}\textbf{Keywords:}\quad Sinai's random walk,  random environment, local time, Brox's diffusion, branching process.\\
\mbox{}\textbf{Mathematics Subject Classification}:  Primary 60J60;
secondary 60G50.
\end{minipage}
\end{center}

\section{ Introduction and Main result\label{s1}}
Let $\{\alpha_{i}\}_{i \in \mathbb{Z}}$ be a sequence of random variables taking values in $(0,1)$. The sequence $\alpha:=\{\alpha_{i}\}_{i \in \mathbb{Z}}$ is called the random environment. For any realization of the environment $\alpha$, we consider the random walk $\{S_{n}\}_{n\geq0}$ defined by $S_{0}=0$,
\begin{equation}
\mathbb{P}(S_{n+1}=i+1 | S_{n}=i,\alpha)=1-\mathbb{P}(S_{n+1}=i+1 | S_{n}=i,\alpha)=\alpha_{i}
\end{equation}
for $n\geq0$ and $i\in\mathbb{Z}$. $\mathbb{P}$ is usually called ``annealed probability". The discrete-time process $\{S_{n}\}_{n\geq0}$ defined in Eq. $(1.1)$ is referred to  Random Walk in Random Environment (RWRE, for short). In this paper, we consider the  Sinai's walk if the random environment satisfies the following assumptions,
\begin{equation}\label{si}
\begin{cases}
\{\alpha_{i}\}_{i\in\mathbb{Z}} ~are ~independent~ and~ identically~ distributed~ variables\\
P(\nu<\alpha_{0}<1-\nu)=1, ~for ~a ~given ~constant ~\nu \in (0,1/2)\\
E(\log\frac{1-\alpha_{0}}{\alpha_{0}})=0\\
0<\sigma^{2}:=E(\log\frac{1-\alpha_{0}}{\alpha_{0}})^{2}<\infty
\end{cases}
\end{equation}
 The corresponding continuous-time model $\{X_{t}\}_{t\geq0}$ has been introduced by Brox~(\cite{B}, 1986), which is often referred to as Brox's diffusion process with Brownian potential. Roughly speaking, it is the solution of the equation
\begin{equation}\label{2}
\begin{cases}
dX_{t}=dB(t)-\frac{1}{2}W'(X_{t})dt, \\
X_{0}=0.
\end{cases}
\end{equation}
where $\{B(t)\}_{t\geq0}$ is one- dimensional standard Brownian motion with $B(0)=0$ and is independent to the medium Brownian motion $\{W(x)\}_{x\in\mathbb{R}}$, and
\begin{equation}
\begin{cases}
W(x)=\sigma W_{1}(x);~~x\geq0, \\
W(x)=\sigma W_{2}(-x);~~x\leq0.
\end{cases}
\end{equation}
where $\{W_{1}(x); x\geq0\}$ and $\{W_{2}(x); x\geq0\}$ are independent one dimensional standard Brownian motion with $W_{1}(0)=W_{2}(0)=0$.
$W'$ denotes the formal derivative of $W$ (it is well-known that Brownian motion is almost surely nowhere differentiable). Rigorously speaking, we can refer $X$ as a Feller-diffusion process on $\mathbb{R}$ with generator
\begin{equation*}
\frac{1}{2}e^{W(x)}\frac{d}{dx}\left(e^{-W(x)}\frac{d}{dx}\right).
\end{equation*}
Actually, for each realization of the potential $\{W(x)\}$, the process $X$ can be represented as
\begin{equation}\label{Brox}
X(t)=A^{-1}(B(T^{-1}(t))),~~t\geq0
\end{equation}
where $A$ is the scale function and $T$ is the time-change function defined by
\begin{align}
&A(y)=\int_{0}^{y}e^{W(z)}dz, ~~y\in\mathbb{R}\label{at}\\
&T(t)=\int_{0}^{t}exp\{-2W(A^{-1}(B(s)))\}ds, ~~t\geq0\label{ts}
\end{align}
Here $A^{-1}$(respectively $T^{-1}$) denotes the inverse of $A$ (respectively $T$).

The significant phenomena for both the Sinai's walk (\cite{S82}) and the Brox diffusion (\cite{B}) are the ``slowly" behavior caused by the random environment,
\begin{equation}
\frac{\sigma^2S_n}{(\ln n)^2}\longrightarrow b_{\infty} \ \ \ \ \mbox{and} \ \ \  \frac{X_t}{(\ln t)^2}\longrightarrow b_{\infty}
\end{equation}
in law for the same variable $b_{\infty}$, whose law has been characterized by  Golosov \cite{G86} and Kesten \cite{Ke86} independently. To get a Donsker's type theorem, Seignourel (\cite{PS}, 2000) designated a discrete scheme in which the medium must be changed at each iteration. In addition the environment $\alpha:=\{\alpha_{i}\}_{i \in \mathbb{Z}}$, define a sequences of environments $\alpha^{(m)}:=\{\alpha_{i}^{(m)}\}_{i\in\mathbb{Z}}$, for $m=1,2,\cdots$,
\begin{equation}\label{am}
\alpha_{i}^{(m)}:=\left(1+\left(\frac{1-\alpha_{i}}{\alpha_{i}}\right)^{\frac{1}{\sqrt{m}}}\right)^{-1},
\end{equation}
and $\{S_{n}^{(m)}\}_{n\geq0}$ is a Random Walk associated with the Random Environment $\alpha^{(m)}:=\{\alpha_{i}^{(m)}\}_{i\in\mathbb{Z}}$, i.e.,
\begin{equation}
\mathbb{P}(S^{(m)}_{n+1}=i+1 | S^{(m)}_{n}=i,\alpha)=1-\mathbb{P}(S^{(m)}_{n+1}=i+1 | S^{(m)}_{n}=i,\alpha)=\alpha^{(m)}_{i}.
\end{equation}

\noindent{\bf Theorem A} (Seignourel, \cite{PS}, 2000) {\it  Under condition (\ref{si}), as $m\to\infty$
\begin{equation}
\left\{ \frac{1}{m}S^{(m)}_{[m^2t]}, t\geq 0 \right\}\longrightarrow \{ X_t, t\geq 0\}
\end{equation}
in distribution in $\mathcal{D}[0,\infty)$, where $\{ X_t, t\geq 0\}$ is the Brox's diffusion process with   $\sigma$ $\cdot$ Brownian motion. \qed

}

 We are interested in the scaling limit of the corresponding local times. For fixed $m\geq 1$, define the local time of $|S^{(m)}|$ at position $j$ before the first $n$ steps as following,
\begin{flalign}\label{lt}
L^{(m)}(j; n)&=\#\{0\leq r\leq n:|S^{(m)}_{r}|=j\} \quad \text{
    for } j,n \geq 0,
\end{flalign}
and define the $n^{th}$ excursion time
\begin{equation}\label{tao}
\begin{cases}
\tau^{(m)}_{0}=0,\\
\tau^{(m)}_{k}=\inf\{n>\tau^{(m)}_{k-1}; |S^{(m)}_{n}|=0\}.
\end{cases}
\end{equation}
For $\forall m\in Z^{+},x\geq 0$, define $L^{(m)}(x)$ as follows:
\begin{equation}\label{lts}
L^{(m)}(x)=
\begin{cases}
\frac{L^{(m)}([mx],\tau^{(m)}_{m})}{m}, & \mbox{for $mx\geq1$}, \\
2, & \mbox{for $0\leq mx <1$}.
\end{cases}
\end{equation}
Now we state our main result.
\begin{theorem}\label{main}
Under condition (\ref{si}), as $m\to\infty$,
\begin{equation}\label{main1}
\{L^{(m)}(x),x\geq0\} ~~ \Rightarrow  ~~\{L^{*}_{X}(x,T(\stackrel{\sim}{T})),,x\geq0\}
\end{equation}
in distribution in $\mathcal{D}[0,\infty)$.
$L^{*}_{X}(x,t)=L_{X}(x,t)+L_{X}(-x,t);~~ t\geq0,x\geq 0$, where $L_{X}(x,t)$ is the local time of the Brox's diffusion process $\{ X_t, t\geq 0\}$ with potential $\sigma$ $\cdot$ Brownian motion and  $\stackrel{\sim}{T}:=\inf\{t\geq0; ~l(0,t)>1\}$, where $l(x,t)$ is the local time of $\{B(t)\}_{t\geq0}$. \qed
\end{theorem}

Donsker's Invariance Principle tells us that the  random walk converge to the
Brownian motion by proper space and time scaling. Accordingly, Rogers (\cite{Rogers}, 1984) have proved that the local time of the simple random walk converges to those    of Brownian motion by scaling based on the branching processes within the random walk (Dwass~\cite{Dwass}, 1975); and is recently generalized to the  random walk with bounded jumps by Hong and Yang (\cite{HY14}, 2014)  by using the intrinsic multi-type branching processes within the random walk (\cite{HW}, 2013). From this point of view, Theorem \ref{main} is a counterpart of the result of Rogers (\cite{Rogers}, 1984) with regard to the  Seignourel's Donsker-type Invariance Principle (Seignourel, \cite{PS}, 2000) for the Sinai's random walk.

We will first express the local time in terms of the branching processes within the random walk (Kesten et al \cite{kks75}, see also Dwass~\cite{Dwass}, 1975); and then prove the convergence of  $\{L^{(m)}(x),x\geq0\}$ by applying the results of  Kurtz~(\cite{Kz},1979) on the convergence of the branching processes in the random environment; at last specify the limiting process being the local time of the Brox's diffusion process  as expected.

\section{ Proof of  Theorem \ref{main}  \label{s3}}

\noindent {\it step 1: Local time and the branching processes within the random walk}

\

For every $m\geq1$, let $(S^{(m)}_{n})_{n\geq0}$ be  random walk in random environment $\{\alpha^{(m)}_{i}\}_{i\in\mathbb{Z}}$, with $S_{0}^{(m)}=0$. Recall the definition (\ref{lt}) and (\ref{tao}),  for the local time $L^{(m)}(j;n)$ of $|S^{(m)}|$ at position $j$ and the excursion  time $\tau^{(m)}_{k}$.
Define, for each $r\geq 1$ and $j\geq 0$,
\begin{equation}
Z_r^{(m)}(j)=\sum_{n=\tau^{(m)}_{r-1}}^{\tau^{(m)}_{r}-1}1_{\{|S^{(m)}_{n}|=j, |S^{(m)}_{n+1}|=j+1\}},
\end{equation}
$Z_r^{(m)}(j)$ records the upward steps from position $j$ to $j+1$ within the $r^{th}$ excursion  of $|S^{(m)}|$.
\begin{proposition}\label{lbp}
(1)For each $r\geq 1$, $(Z_r^{(m)}({j}))_{j\geq0}$ is a branching process in random environment with $Z_r^{(m)}({0})=1$ and the offspring branching mechanism
\begin{equation}\label{be1}
\mathbb{P}\left( {Z_r^{(m)}({j+1})=k}\mid Z_r^{(m)}({j})=1; ~\alpha^{(m)}_{i},i\in\mathbb{Z}\right)=\left( \alpha^{(m)}_{j+1} \right)^{k}\left(1-\alpha^{(m)}_{j+1}\right).
\end{equation}
i.e.,
\begin{equation}\label{be}
\mathbb{E}\left(s^{Z_r^{(m)}({j+1})}\mid Z_r^{(m)}({0}),Z_r^{(m)}({1}),\cdots ,Z_r^{(m)}({j}); ~\alpha^{(m)}_{i},i\in\mathbb{Z}\right)=\left(\frac{1-\alpha^{(m)}_{j+1}}{1-\alpha^{(m)}_{j+1}\cdot s}\right)^{Z_r^{(m)}({j})}.
\end{equation}
(2)Within the first excursion, the local time at position $j\geq 1$ can be expressed in terms of  the branching process in random environment $(Z^{(m)}_{j})_{j\geq0}$,
\begin{equation}\label{ex1}
 L^{(m)}(j; \tau^{(m)}_{1}) = Z_1^{(m)}({j-1})+Z_1^{(m)}({j}).
\end{equation}
(3) For any positive integral $N$,
\begin{equation}\label{exn}
 L^{(m)}(j; \tau^{(m)}_{N}) = \sum_{r=1}^{N} \left(L^{(m)}(j; \tau^{(m)}_{r})-L^{(m)}(j; \tau^{(m)}_{r-1})\right):=\sum_{r=1}^{N}\left(Z_r^{(m)}({j-1})+Z_r^{(m)}({j})\right),
\end{equation}
where  for any fixed environment, $\eta_r^{(m)}({j}):=\left(Z_r^{(m)}({j-1})+Z_r^{(m)}({j})\right)$, $r=1,2,\cdots$ are i.i.d. random variables, distributed as $L^{(m)}(j; \tau^{(m)}_{1})$.\qed
\end{proposition}
Part (1)   is from Kesten et al (\cite{kks75}, 1975), see also Dwass (\cite{Dwass}, 1975) for the simple random walk; Part (2) and (3) is similar as    {Theorem 1}  of Rogers~(\cite{Rogers}, 1984).

Proposition \ref{lbp} tells us that for each $r\geq 1$,  $(Z_r^{(m)}({j}))_{j\geq0}$ is a branching process in random environment with $Z_r^{(m)}({0})=1$; and they are independent of each other.  If we write
\begin{equation}\label{U}
U_N^{(m)}(j)=\sum_{r=1}^{N}Z_{r}^{(m)}(j),
\end{equation}
 $\{U_N^{(m)}(j)\}_{j\geq0}$  is a branching process in random environment with $U_N^{(m)}({0})=N$ and the branching mechanism given in (\ref{be1}).  From (\ref{exn}) and (\ref{U}), for any positive integral $N$, we have
\begin{equation}\label{exn1}
 L^{(m)}(j; \tau^{(m)}_{N}) =\sum_{r=1}^{N}\left(Z_r^{(m)}({j-1})+Z_r^{(m)}({j})\right)=U_N^{(m)}(j-1)+U_N^{(m)}(j).
\end{equation}
As a consequence, for $\forall m\in Z^{+},x\geq 0$,  $L^{(m)}(x)$, the scaling of the local time  (see (\ref{lts})) can be written as,
\begin{flalign}
L^{(m)}(x)&:=\frac{L^{(m)}([mx],\tau^m_{m})}{m}\nonumber\\
&=\frac{U_{m}^{(m)}([mx])}{m}+\frac{U_{m}^{(m)}([m(x-m^{-1})])}{m}.
\end{flalign}
Define
\begin{equation}\label{X}
X_{m}(x)=\frac{U_{m}^{(m)}([mx])}{m}
\end{equation}
then $L^{(m)}(x)=X_{m}(x)+X_{m}(x-m^{-1})$, so (\ref{main1}) in Theorem 1.1 will follow if we can prove
\begin{equation}\label{d}
2X_{m}(\cdot)~~\Rightarrow~~L^{*}_{X}(\cdot,T(\stackrel{\sim}{T}))
\end{equation}

\

\noindent {\it step 2: convergence of  $\{X_{m}(x),x\geq0\}$}

\

Now we will consider the  convergence of  $\{X_{m}(x),x\geq0\}$ in distribution in $\mathcal{D}[0,\infty)$,
 based on Theorem (2.13) of Kurtz~(\cite{Kz},1979).  Firstly from (\ref{U}), we know that $\{U_m^{(m)}(k)\}_{k\geq0}$  is a branching process in random environment with $U_m^{(m)}({0})=m$, i.e., for $k\geq 1$,
 \begin{flalign}\label{U1}
&U_m^{(m)}(k)=\sum_{i=1}^{U_m^{(m)}(k-1)}\xi_{i,k}^{(m)}
\end{flalign}
where $\{\xi_{i,k}^{(m)}\}_{i\geq 1}$ are i.i.d random variables whose distribution is given in (\ref{be1}).
With (\ref{X}) in mind, we will check the conditions of Theorem (2.13) of Kurtz~(\cite{Kz},1979) directly. To this end, recall the notations in (\cite{Kz}), we can calculate the moments of $\xi_{i,k}^{(m)}$ from the branching mechanism (\ref{be1}) or (\ref{be}),
\begin{flalign}\label{1m}
&\lambda_{k}^{(m)}=\mathbb{E}(\xi_{i,k}^{(m)}|~\{\alpha^{(m)}_{j}\}_{j\in\mathbb{Z}})= \frac{ \alpha^{(m)}_{k}}{1-\alpha^m_{k}}:=\rho^{(m)}_{k},
\end{flalign}
note that from (\ref{am}), when $m=1$ we have $\{\alpha^{(m)}_{j}\}_{j\in\mathbb{Z}}=\{\alpha_{j}\}_{j\in\mathbb{Z}}$, and we write $ \rho_{k}=\frac{ \alpha_{k}}{1-\alpha_{k}}$;  with this notation,  (\ref{am}) can be rewritten as
\begin{equation}\label{am1}
\alpha_{i}^{(m)}:=\left(1+\left(\rho_{i}\right)^{-\frac{1}{\sqrt{m}}}\right)^{-1},
\end{equation}
and as a consequence the first moment (\ref{1m}) equals to
\begin{flalign}\label{1m1}
&\lambda_{k}^{(m)}=\mathbb{E}(\xi_{i,k}^{(m)}|~\{\alpha^{(m)}_{j}\}_{j\in\mathbb{Z}})= \frac{ \alpha^{(m)}_{k}}{1-\alpha^m_{k}}:=\rho^{(m)}_{k}=(\rho_{k})^{\frac{1}{\sqrt{m}}},
\end{flalign}
For the second moment
\begin{flalign}\label{2m}
&a_{k}^{(m)}=\mathbb{E}\left((\xi_{i,k}^{(m)}/\lambda_{k}^{(m)}-1)^{2}|~\{\alpha^{(m)}_{i}\}_{i\in\mathbb{Z}}\right)
=\left(\alpha_{i}^{(m)}\right)^{-1}=\left(1+\left(\rho_{i}\right)^{-\frac{1}{\sqrt{m}}}\right),
\end{flalign}
and the third moment
\begin{flalign}\label{3m}
Y_{k}^{(m)}&=\mathbb{E}(|\xi_{i,k}^{(m)}/\lambda_{k}^{(m)}-1|^{3}|~\{\alpha^{(m)}_{i}\}_{i\in\mathbb{Z}})\nonumber\\
&\leq{\mathbb{E}((\xi_{i,k}^{(m)}/\lambda_{k}^{(m)}+1)^{3}|~\{\alpha^{(m)}_{i}\}_{i\in\mathbb{Z}})}\nonumber\\
&\leq {2^{2}[\mathbb{E}(\xi_{i,k}^{(m)}/\lambda_{k}^{(m)})^{3}|~\{\alpha^{(m)}_{i}\}_{i\in\mathbb{Z}})+1]}\nonumber\\
&\leq4{(7+6(\rho_{k}^{-1})^{\frac{1}{\sqrt{m}}}+(\rho_{k}^{-2})^{\frac{1}{\sqrt{m}}})},
\end{flalign}
at the last step it's easy to calculate that $
\mathbb{E}(\xi_{i,k}^{(m)}/\lambda_{k}^{(m)})^{3}|~\{\alpha^{(m)}_{i}\}_{i\in\mathbb{Z}})=6+6/\rho_{k}^{(m)}+1/(\rho_{k}^{(m)})^{2}$.
Now we can check the three conditions of Theorem (2.13) of Kurtz~(\cite{Kz},1979).

\noindent{\it (1) condition 1:}  From (\ref{1m1}), we have
\begin{flalign}\label{cd1}
 M_{m}(x)&:=\sum_{k=1}^{[mx]}\log\lambda_{k}^{(m)}=\frac{1}{\sqrt{m}}\sum_{k=1}^{[mx]}\log\rho_{k},
\end{flalign}
then as a process $\{M_{m}(x)\}_{x\geq0}\Rightarrow \{M(x)\}_{x\geq0}$ when $m\rightarrow \infty$ by Donsker's invariance principle with regard of the condition (\ref{si}), where $\{M(t)\}_{t\geq0}\stackrel{d}{=}\{\sigma B(t)\}_{t\geq0}$, $\{B(t)\}_{t\geq0}$ is a standard Brownion Motion.

\noindent{\it (2) condition 2:} From (\ref{2m}), for any $x\geq 0$,
\begin{flalign}\label{cd2}
 A_{m}(x):&=\frac{1}{m}\sum_{k=1}^{[mx]}a_{k}^{(m)}
 =\frac{1}{m}\sum_{k=1}^{[mx]}\left(1+\left(\rho_{i}\right)^{-\frac{1}{\sqrt{m}}}\right)\nonumber\\
 &\longrightarrow 2x
\end{flalign}
$P$-a.s. as $m\to\infty$ by the elliptic condition for the environment in (\ref{si}).

\noindent{\it (3) condition 3:} From (\ref{3m}), for any $x\geq 0$,
\begin{flalign}\label{cd3}
G_{m}(x)&:=\frac{1}{m^{3/2}}\sum_{k=1}^{[mx]}{Y_{k}^{m}}\leq\frac{4}{m^{3/2}}
\sum_{k=1}^{[mx]}{(7+6(\rho_{k}^{-1})^{\frac{1}{\sqrt{m}}}+(\rho_{k}^{-2})^{\frac{1}{\sqrt{m}}})}\nonumber\\
 &\longrightarrow 0
\end{flalign}
$P$-a.s. as $m\to\infty$ by the elliptic condition for the environment in (\ref{si}).

Combining the above three conditions {\it{(1)-(3)}}, by Theorem (2.13) of Kurtz~(\cite{Kz},1979), we have proved the following
\begin{theorem}\label{xl} As $m\to\infty$,  we have
\begin{flalign}\label{XH}
\{X_{m}(x)\}_{x\geq0}~~\Rightarrow ~~\{H(x)\}_{x\geq0}:=\{V(\beta(x))e^{M(x)}\}_{x\geq0},
\end{flalign}
where $\{V(x)\}_{x\geq0}$ is a diffusion process with generator $\frac{1}{2}yf^{\prime\prime}(y)$, which is independent of $\{M(x)\}_{x\geq0}$; and  $\beta(x)=\tau^{-1}(x)$,   where  $\tau(x)$ is determined by
\begin{flalign}\label{tau}
x &=\int_{0}^{\tau(x)}{(2s)^{\prime}e^{-M(s)}ds}=\int_{0}^{\tau(x)}{2e^{-M(s)}ds}.
\end{flalign}\qed
\end{theorem}

\noindent{\bf Remark 2.1} {\it From
Chapter 9 Theorem 3.1 in Kurtz~(\cite{Kurtz},2005) , see also Kurtz~(\cite{Kz},1979),  we know that the
Markov process $H$ is the unique solution of
\begin{equation*}
H(t)=e^{M(t)}[1+B(\int_{0}^{t}{e^{-2M(s)}H(s)d(2s)})]
\end{equation*}
and the generator is
\begin{equation*}
Af(z)=\frac{1}{2}\sigma^{2}zf^{\prime}(z)+(z+\frac{1}{2}\sigma^{2}z^{2})f^{\prime\prime}(z)
\end{equation*}
for every $f\in C_{c}(\mathbb{R})$. \qed
}

\

\noindent{\bf Remark 2.2} {\it  $\{V(x)\}_{x\geq0}$ is a diffusion process with generator $\frac{1}{2}yf^{\prime\prime}(y)$, i.e., is the solution of the stochastic differential equation with $V(0)=1$
\begin{flalign}\label{veq}
&V(x)=1+ \int_{0}^{x}(V(s)^{+})^{\frac{1}{2}}dB_{s},
\end{flalign}
 for $a>0$, it is easy to verify that  $\{V(x)\}_{x\geq0}$ has the same generator with $ \{aV\left(\frac{x}{a}\right)\}_{x\geq0}$. As a consequence, \{$\stackrel{\sim}{V}(x)\}_{x\geq0}:=\{V(2x)\}_{x\geq0}$ satisfies the stochastic differential equation with $\stackrel{\sim}{V}(0)=1$
\begin{flalign}\label{veq1}
&\stackrel{\sim}{V}(x)=1+ \sqrt{2}\int_{0}^{x}(\stackrel{\sim}{V}(s)^{+})^{\frac{1}{2}}dB_{s},
\end{flalign}
Then by Ray-Knight theorem we know that
\begin{equation}\label{lv}
\{l(x; {\stackrel{\sim}{T}})\}_{x\geq0}\stackrel{d}{=}\{\stackrel{\sim}{V}(x)\}_{x\geq0}
\end{equation}
where $\stackrel{\sim}{T}=\inf\{t:l(0,t)>1\}$, $l(x; t)$ is the local time at position $x$ of the standard Brownian motion. \qed

}

\noindent{\bf Remark 2.3} {\it From (\ref{tau}), we know  $\beta(x)=\tau^{-1}(x)$  can be expressed as
\begin{flalign}\label{tau1}
\beta(x)&=2\int_{0}^{x}{e^{-M(s)}ds}\stackrel{d}{=}2 A(x),
\end{flalign}
where $A(x)$ is given in (\ref{at}) because $\{M(t)\}_{t\geq0}\stackrel{d}{=}\{\sigma B(t)\}_{t\geq0}$, $\{B(t)\}_{t\geq0}$ is a standard Brownion Motion. So the limiting process $\{H(x)\}_{x\geq0}$ in (\ref{XH}) has the following
\begin{flalign}\label{XH1}
\{H(x)\}_{x\geq0}:&=\{V(\beta(x))e^{M(x)}\}_{x\geq0}\stackrel{d}{=}
\{V(2A(x))e^{-M(x)}\}_{x\geq0}\stackrel{d}{=}\{\stackrel{\sim}{V}(A(x))e^{-M(x)}\}_{x\geq0}\nonumber\\
& \stackrel{d}{=}\{l(A(x); \stackrel{\sim}{T})e^{-M(x)}\}_{x\geq0},
\end{flalign}
the last step is by (\ref{lv}) where $\stackrel{\sim}{T}=\inf\{t:l(0,t)>1\}$.
\qed
}

\

\noindent {\it step 3: local time of Brox diffusion}

\

Denote $L_{X}(x; t)$ as the   local time of Brox diffusion $\{X(t)\}_{t\geq0}$ with potential $\{M(x)\}$.  To find $L_{X}(x,t)$,  according to the definition of the local time,  it is enough for any bounded Borel function $f$,
\begin{equation*}
\int_{0}^{t}{f(X(s))ds}=\int_{-\infty}^{+\infty}{f(x)L_{X}(x; t)dx}.
\end{equation*}
 Recall $l(x; t)$ is the local time process of standard Brownion motion, by (\ref{Brox}), (\ref{at}) and (\ref{ts}) we have (see for example, \cite{HS})
\begin{flalign*}
\int_{0}^{t}{f(X(s))ds}&=\int_{0}^{t}{f(A^{-1}(B(T^{-1}(s))))ds}\\
&=\int_{0}^{T^{-1}(t)}{f(A^{-1}(B(u)))\exp[-2M(A^{-1}(B(u)))]du}\\
&=\int_{-\infty}^{+\infty}{f(A^{-1}(y))\exp(-2M(A^{-1}(y)))l(y; T^{-1}(t))dy}\\
&=\int_{-\infty}^{+\infty}{f(x)\exp(-2M(x))l(A(x),T^{-1}(t))e^{M(x)}dx}\\
&=\int_{-\infty}^{+\infty}{f(x)e^{-M(x)}l(A(x),T^{-1}(t))dx}.
\end{flalign*}
Hence, we get
\begin{equation}\label{bl}
L_{X}(x; t)=l(A(x),T^{-1}(t))e^{-M(x)},
\end{equation}
where $A$ and $T$ has defined in~\eqref{at},~\eqref{ts}.\\

\

\noindent {\it step 4: conclusion}

\

By the equation~\eqref{bl}, if we let $t=T(\stackrel{\sim}{T})$, then,
\begin{equation}\label{c}
L_{X}(x,T(\stackrel{\sim}{T}))=l(A(x),\stackrel{\sim}{T})e^{-M(x)}.
\end{equation}
Combining  (\ref{XH1}) we have
\begin{flalign}\label{XH2}
\{H(x)\}_{x\geq0}:&=\{V(\beta(x))e^{M(x)}\}_{x\geq0}  \stackrel{d}{=}\{l(A(x); \stackrel{\sim}{T})e^{-M(x)}\}_{x\geq0}
 \stackrel{d}{=}L_{X}(x,T(\stackrel{\sim}{T})),
\end{flalign}
where $\stackrel{\sim}{T}=\inf\{t:l(0,t)>1\}$. Thus (\ref{XH}) in Theorem \ref{xl} is actually
\begin{flalign}\label{XH3}
\{X_{m}(x)\}_{x\geq0}~~\Rightarrow ~~\{L_{X}(x,T(\stackrel{\sim}{T}))\}_{x\geq0}.
\end{flalign}
To finish the proof of (\ref{d}), we note that by the property of the local time of Brownian motion,
\begin{equation*}
\{2 l(x,\stackrel{\sim}{T})\}_{x\geq0}\stackrel{d}
{=}\{l(x,\stackrel{\sim}{T})\}_{x\geq0}+\{l(-x,\stackrel{\sim}{T})\}_{x\geq0},
\end{equation*}
we have
\begin{flalign*}
\{L^{*}_{X}(x,T(\stackrel{\sim}{T}))\}_{x\geq0}&\stackrel{d}
{:=}\{L_{X}(x,T(\stackrel{\sim}{T}))+L_{X}(-x,T(\stackrel{\sim}{T}))\}_{x\geq0}\\
&=\{l(A(x),\stackrel{\sim}{T})e^{-M(x)}+l(A(-x),\stackrel{\sim}{T})e^{-M(-x)}\}_{x\geq0}\\
&=\{l(A(x),\stackrel{\sim}{T})e^{-M(x)}+l(-\int_{0}^{x}{e^{M(-y)}}dy,\stackrel{\sim}{T})e^{-M(-x)}\}_{x\geq0}\\
&\stackrel{d}{=}\{l(A(x),\stackrel{\sim}{T})e^{-M(x)}+l(-\int_{0}^{x}{e^{M(y)}}dy,\stackrel{\sim}{T})e^{-M(x)}\}_{x\geq0}\\
&=\{l(A(x),\stackrel{\sim}{T})e^{-M(x)}+l(-A(x),\stackrel{\sim}{T})e^{-M(x)}\}\\
&\stackrel{d}{=}\{2l(A(x),\stackrel{\sim}{T})e^{-M(x)}\}_{x\geq0}\\
&\stackrel{d}{=}\{2V(\beta(x))e^{M(x)}\}_{x\geq0}.
\end{flalign*}
which complete the proof of Theorem 1.1. \qed

\end {document}